\documentclass[preprint,12pt]{elsarticle}

\usepackage{geometry}                
\usepackage{graphicx}
\usepackage{amssymb,amsmath,amsthm}
\usepackage{bbm}
\usepackage{color}

\newtheorem{theorem}{Theorem}[section]

\newtheorem{lemma}[theorem]{Lemma}

\newtheorem{proposition}[theorem]{Proposition}

\newtheorem{definition}[theorem]{Definition}

\newtheorem{example}[theorem]{Example}

\newtheorem{remark}[theorem]{Remark}

\newproof{pf}{Proof}
\newproof{pot}{Proof of Theorem \ref{thm2}}
\renewcommand{\qedsymbol}{$\blacksquare$}

\newcommand{\N}{\mathbb{N}}

\newcommand{\R}{\mathbb{R}}

\newcommand{\ep}{\epsilon }
\newcommand{\la}{\lambda }

\newcommand{\de}{\delta }

\newcommand{\si}{\sigma }

\newcommand{\ga}{\gamma }
\newcommand{\Ga}{\Gamma }

\newcommand{\dess}{d_{\text ess}}

\newcommand{\length}{\operatorname{length}}

\newcommand{\Mod}{\operatorname{Mod}}

\newcommand{\adm}{\operatorname{Adm}}
\newcommand{\wadm}{\operatorname{w-Adm}}

\newcommand{\defeq}{\mathrel{\mathop:}=}

\newcommand{\ones}{\mathbbm{1}}

\title{Infinity modulus and the essential metric\tnoteref{t1}}
\tnotetext[t1]{The first and third author were supported by NSF grant DMS-\#1515810. The fourth author was supported by NSF grant 
DMS-\#1500440.}

\author[ksu]{Nathan Albin}
\ead{albin@math.ksu.edu}

\author[ksu]{Jared Hoppis}
\ead{jhoppis@math.ksu.edu}

\author[ksu]{Pietro Poggi-Corradini\corref{cor1}\fnref{fn1}}
\ead{pietro@math.ksu.edu}
\fntext[fn1]{The third named author thanks Marshall Williams for an interesting discussion}
\cortext[cor1]{Corresponding author}

\author[cin]{Nageswari Shanmugalingam}

\ead{shanmun@uc.edu}

\address[ksu]{Department of Mathematics, Cardwell Hall, Kansas State University,
Manhattan, KS 66506, USA}
\address[cin]{Department of Mathematical Sciences, P.O. Box 210025, University of Cincinnati,
Cincinnati, OH~45221-0025, USA}

\begin{document}

\begin{abstract}
We study $\infty$-modulus on general metric spaces and establish its relation to shortest lengths of paths. This connection was already known for modulus on graphs, but the formulation in metric measure spaces requires more attention to exceptional families. We use this to define a metric that we call the essential metric, and show how this recovers a metric that had already been  advanced in the literature by De Cecco and Palmieri.  
\end{abstract}

                     
\maketitle
\baselineskip=18pt

\section{Introduction}

The $p$-modulus of a family of curves  is a way to quantify the richness of such a family. This began as an important tool in function theory, because of its conformal invariance in the $p=2$ case~\cite{ahlfors1973,beurling1989}. Later, $p$-modulus was defined and used successfully on $n$-dimensional Euclidean spaces, and general metric spaces as well~\cite{heinonen2001,hkst2015}. The case $p=\infty$ has been studied on general metric spaces by the fourth author and her collaborators~\cite{Durandcartagena-Jaramillo-Shanmugalingam:CVGMT2016,durandcartagena-jaramillo:jmaa2010}, and on graphs and networks by the first and third authors and their collaborators~\cite{abppcw:ecgd2015}. In particular, on graphs, $\infty$-modulus was found to be connected to the shortest path graph distance~\cite[Theorem 4.1]{abppcw:ecgd2015}. This fact gave the impetus for the present paper. Here, we develop the theory of $\infty$-modulus on general metric space from first principles. The presentation is mostly self-contained, so that anyone with some background in analysis can follow easily.
We pursue the connection between $\infty$-modulus and shortest path length in this more abstract setting and we define a new notion of {\it essential length} of a family of curves. This then gives rise to a new metric that we call the {\it essential metric}. In the process, we revisit and reinterpret some of the existing results from the literature in this new light.

The essential metric is one way to measure the effective shortest
path between two points. In the case of a graph, it is indeed associated with the shortest path between two points. In the setting of metric
spaces the trajectory of a single curve (or countably many curves) might have measure zero; this does not happen in the case of 
a graph or network where the underlying space is locally one-dimensional. Thus, in the metric setting the essential metric does not measure the absolute
shortest path between the two points. Instead, it can ignore a subfamily  of ``shortest'' paths, if that  collection is negligeable (that is, has zero modulus). For example, the Sierpinski gasket is a quasiconvex metric
space when viewed as a metric subspace of the Euclidean two-dimensional space and equipped with the 
$\log(3)/\log(2)$-dimensional Hausdorff measure, but the essential distance between distinct pairs of
points there is infinite. Hence, in general, it is of interest to know which pairs of points have finite essential distance, and which pairs of points do not.

This paper is structured as follows. In
Section \ref{sec:basic},   we recall the
the basic tools necessary to compute the length of a curve in a metric space.
Then in Section \ref{sec:infty-mod}, we give two equivalent definitions of $\infty$-modulus on metric measure spaces,
characterize the notion of $\infty$-exceptional families, 
introduce the notion of essential length of a family of curves, and 
show how it relates to $\infty$-modulus.
Finally, in Section \ref{sec:ess-metric}, we define the essential metric and show that it coincides with a differently defined metric  that was first introduced by De Cecco and Palmieri.

\section{Basic tools in metric spaces}\label{sec:basic}
\subsection{Length in metric spaces}

We follow Chapter 5 in \cite{hkst2015}. The material in this section can be found in other sources as well, but we collect as much as possible here, for the reader's convenience, so as not to require a lot of background.

A {\it metric space} $(X,d)$ is a set $X$ equipped with a {\it metric} $d$. That is, $d$ is a function $d: X\times X\rightarrow\R$ satisfying non-degeneracy: $d(x,y)\ge 0$ for all $x,y\in X$, and $d(x,y)=0$ if and only if $x=y$ ; symmetry: $d(x,y)=d(y,x)$ for all $x,y\in X$; and finally, the triangle inequality:
\[
d(x,y)\le d(x,z)+d(z,y)\qquad\forall x,y,z\in X.
\]

A {\it path} $\ga$ in $X$ is a continuous map $\ga:[a,b]\longrightarrow (X,d )$ for some
compact interval $[a,b]\subset\mathbb{R}$.  Its {\it length} is the total variation:
\begin{equation}\label{eq:length}
\length(\ga)\defeq \sup_{a=t_{0}\leq t_{1}\dots\leq t_{N}=b}\sum_{k=0}^{N-1}d(\ga(t_{k}),\ga(t_{k+1})),
\end{equation}
where the supremum is taken over all possible partitions with $N$
arbitrary. 

The path $\ga$ is called {\it rectifiable} if $\length (\ga)<\infty$. In this case, we can define the {\it length function} $s_{\ga}:[a,b]\longrightarrow[0,\length(\ga)]$ as
\begin{equation}
s_{\ga}(t)\defeq \length\left(\ga|_{[a,t]}\right).
\end{equation}
Clearly, $s_{\ga}$ is increasing (in the weak sense). Also it can be checked using (\ref{eq:length}) that for any $a\leq t<s\leq b$:
\begin{equation}\label{eq:distlength}
d(\ga(t),\ga(s))\leq \length\left(\ga|_{[t,s]}\right)=s_{\ga}(s)-s_{\ga}(t).
\end{equation}

\begin{lemma}\label{lem:sgacont}
If $\ga$ is a rectifiable path in $X$, then $s_{\ga}$ is continuous.
\end{lemma}

\begin{pf}
Suppose $s_\ga$ is not left continuous. Then there is $a<t_{0}\le b$
and $\delta>0$ such that 
\begin{equation}\label{eq:jump}
s_{\ga}(t_{0})>s_{\ga}(t)+\delta \qquad \text{for all $a<t<t_{0}.$}
\end{equation}
Since $\ga$ is continuous at $t_{0}$, there is $a<t_{1}<t_{0}$ so that 
\begin{equation}\label{eq:cont}
d(\ga(t),\ga(t_{0}))<\delta/2\qquad \text{for every $t_{1}<t<t_{0}$.}
\end{equation}
By (\ref{eq:jump}), $\length\left(\ga|_{[t_{1},t_{0}]}\right)=s_{\ga}(t_{0})-s_{\ga}(t_{1})>\de$. So, there is a 
partition $t_{1}=s_{0}<s_{1}<\cdots<s_{N-1}<s_{N}=t_{0}$ of $[t_1,t_0]$, so that
\begin{equation}
\sum_{j=0}^{N-1}d(\ga(s_{j}),\ga(s_{j+1}))>\delta.
\end{equation}
Using (\ref{eq:cont}) with $t=s_{N-1}$, we get that
\begin{equation}
\sum_{j=0}^{N-2}d(\ga(s_{j}),\ga(s_{j+1}))>\delta/2.
\end{equation}
Let $t_{2}\defeq s_{N-1}$ and repeat the steps above with $t_{2}$ in place of $t_{1}$, noting that both (\ref{eq:jump}) and (\ref{eq:cont}) remain valid.
This process can be iterated indefinitely and it gives rise to a sequence $t_{k}<t_0$, such that $\length\left(\ga|_{[t_{1},t_{k}]}\right)\rightarrow\infty$. This contradicts the rectifiability of $\ga$.

The proof that $s_\ga$ is right continuous is similar.
\hfill\qedsymbol
\end{pf}

The length function $s_\ga$ is continuous and increasing, but not necessarily strictly increasing. Still, we can define a  right-inverse 
$s_\ga^{-1}:[0,\length\left(\ga\right)]\to[a,b]$
as follows:
\begin{equation}\label{eq:sgainv}
s_\ga^{-1}(t)=\max\{s: s_\ga(s)=t\}\qquad\forall t\in [0,\length(\ga)]
\end{equation}
Then $s_\ga^{-1}$ is increasing, right-continuous, and $s_\ga(s_\ga^{-1}(t))=t$.
\begin{definition}
The {\rm arc-length parametrization} of a rectifiable path $\ga:[0,1]\rightarrow X$ is the curve $\ga_s:[0,\length(\ga)]\rightarrow X$ defined by
\[
\ga_s(t):=\ga(s_\ga^{-1}(t)).
\] 
\end{definition}
In particular, $\ga(u)=\ga_s(s_\ga(u))$, and
\begin{equation}\label{eq:arclength}
\length\left(\ga_s|_{[t,u]}\right)=\length\left(\ga|_{[s_\ga^{-1}(t),s_\ga^{-1}(u)]}\right)=s_{\ga}(s_\ga^{-1}(t))-s_{\ga}(s_\ga^{-1}(u))=t-u.
\end{equation}

\begin{definition}\label{def:abscont}
A path
$\ga:[0,1]\rightarrow X$ is {\rm absolutely continuous} if for all $\epsilon>0$ there exists $\de=\de(\ep)>0$ such that whenever $\{(a_i,b_i)\}_{i=1}^N$ are disjoint intervals in $[0,1]$:
\[
\sum_{i=1}^N |b_i-a_i|<\de\Longrightarrow \sum_{i=1}^N d(\ga(a_i),\ga(b_i))<\ep.
\]
\end{definition}

\begin{proposition}\label{prop:abs-cont}
Suppose $\ga:[0,1]\rightarrow X$ is rectifiable. Then $\ga$ is absolutely continuous if and only if $s_\ga$ is absolutely continuous. 
\end{proposition}
\begin{pf}
If $s_\ga$ is absolutely continuous, then by (\ref{eq:distlength}), the same is true of $\ga$.

Conversely, assume that $\ga$ is absolutely continuous.  Given $\ep>0$, find $\de=\de(\ep)>0$ as in Definition \ref{def:abscont}. Let $\{(a_i,b_i)\}_{i=1}^N$ be disjoint intervals in $[0,1]$ with
\[
\sum_{i=1}^N |b_i-a_i|<\de.
\]
Then, as we have seen, $s_\ga(b_i)-s_\ga(a_i)=\length\left(\ga\mid_{[a_i,b_i]}\right)<\infty$.
By (\ref{eq:length}), there are disjoint intervals $\{(a_i^j,b_i^j)\}_{j=1}^{N_i}$ contained in $(a_i,b_i)$ such that
\[
\sum_{j=1}^{N_i}d(\ga(a_i^j),\ga(b_i^j))\geq s_\ga(b_i)-s_\ga(a_i)-\frac{\ep}{N}.
\]
By absolute continuity of $\ga$, since the disjoint intervals $\{(a_i^j,b_i^j)\}_{i,j}$ also have length adding up to less than $\de$, we get
\[
\sum_{i=1}^N s_\ga(b_i)-s_\ga(a_i)\leq \sum_{i=1}^N\sum_{j=1}^{N_i}d(\ga(a_i^j),\ga(b_i^j)) +\ep\leq 2\ep.
\]
\hfill\qedsymbol
\end{pf}

Next we recall Proposition 5.1.8 of \cite{hkst2015}.

\begin{proposition}[\cite{hkst2015}]\label{prop:hkst}
Let $\ga:[0,1]\rightarrow X$ be a compact rectifiable path. Then, its arc-length parametrization $\ga_s$ is absolutely continuous. Indeed, $\ga_s$ is $1$-Lipschitz and
\[
\lim_{u\rightarrow t, u\neq t}\frac{d(\ga_s(t),\ga_s(u))}{|t-u|}=1\qquad\text{for a.e. $t$.}
\]
\end{proposition}

\begin{definition}
Suppose $\ga:[0,1]\rightarrow X$ is rectifiable and $\rho:X\rightarrow[0,\infty]$ is Borel. Then the {\rm line integral} of $\rho$ along $\ga$ is
\begin{equation}\label{eq:line-integral}
\int_\ga \rho ds := \int_0^{\length(\ga)}\rho(\ga_s(t))dt.
\end{equation}
Also, if $F\subset X$ is a Borel set, we say that $\ga$ {\rm has positive length in} $F$ if
\[
\int_{\ga\cap F} ds := \int_\ga \ones_F\, ds =  \int_0^{\length(\ga)} \ones_{\ga_s^{-1}(F)}(t)dt=m_1\left(\ga_s^{-1}(F)\right) >0,
\]
where $m_1$ is the Lebesgue measure on $\R$.
We write $\Ga_F^\ell$ for the family of all curves that have positive length in $F$.
\end{definition}

The key observation for (\ref{eq:line-integral}) is that the composition of the Borel function $\rho$ and the continuous function $\ga_s$, is a measurable function.

In this paper, a {\it curve} will denote a non-constant path, defined on a possibly infinite interval $[a,b]$, that is locally rectifiable, meaning that every $t\in(a,b)$ has a neighborhood where $\gamma$ is rectifiable.  Unless otherwise stated,  all the families that will be considered will be families of such curves.
\begin{definition}
Suppose $\ga:[0,1]\rightarrow X$ is a curve and  $F\subset X$ is a Borel set.
We say that  $\ga$ {\rm spends positive time} in $F$  if 
\[
\int_0^{1} \ones_{\ga^{-1}(F)}(t)dt=m_1\left(\ga^{-1}(F)\right) >0.
\]
We write $\Ga_F^\tau$ for the family of all curves that spend positive time in $F$.
\end{definition}
The two concepts of having positive length in $F$ and spending positive time in $F$ are in general unrelated.
For instance,  suppose $\ga$ is a curve traveling from left to right on $\R$ at constant speed, and suppose $\ga$ stops at the origin for one unit of time. Then,
if $F=\{0\}$ is the singleton containing the origin, $\ga$ spends positive time in $F$, but $\ga$ does not have positive length in $F$.
Conversely, consider the curve
\[
\ga(t)=(t,C(t)),\qquad t\in[0,1],
\]
where $C(t)$ is the usual Cantor step-function. Let $C$ be the middle-third Cantor set, $D$ be the dyadic rationals in $[0,1]$, and $F = [0,1] \times([0,1] \backslash D )$. Then $\ga^{-1}(F)=C$, and so $m_1\left(\ga^{-1}(F)\right) =0.$ Intuitively, the curve $\ga$ has infinite speed on $C$, and therefore spends zero time on $\ga(C)$.
Now, letting $E =[0,1] \times D $, we have 
\[
\sqrt{2} \leq \length(\ga) = m_1\left(\ga_s^{-1}([0,1] \times [0,1])\right)= m_1\left(\ga_s^{-1}(F)\right) + m_1\left(\ga_s^{-1}(E)\right)
\]
A simple computation shows $m_1\left(\ga_s^{-1}(E)\right)= 1$, and hence, $m_1\left(\ga_s^{-1}(F)\right) \geq \sqrt{2} - 1 > 0$. Therefore, $\ga$ has positive length in $F$. Note also that $F$ can be  intersected with $\ga([0,1])$ so as to get an example where $F$ has area measure zero.

On the other hand, if $\ga$ is absolutely continuous, then for every Borel set $F\subset X$ we have
\[
m_1\left(\ga^{-1}(F)\right) =0\Longrightarrow m_1\left(\ga_s^{-1}(F)\right) =0,
\]
meaning that, in this case, if $\ga$ has positive length in $F$, then it also spends positive time in $F$.
To see this, let $f = s_{\ga}^{-1}$ and $g = \ga$. Then we have $f^{-1}(A) \subset s_{\ga}(A)$, $\forall A \subset [0,1]$ and also,
\[
 \ga_{s}^{-1} (F) = \left( g \circ f \right)^{-1}(F) = f^{-1} \left( g^{-1}(F) \right).
\]
By Proposition \ref{prop:abs-cont}, $s_{\ga}$ is absolutely continuous, and hence,
\[
m_1\left(\ga^{-1}(F)\right) =0 \Rightarrow m_1\left( s_{\ga}\left( \ga^{-1}(F) \right)\right) =0 \Rightarrow m_1\left(\ga_{s}^{-1}(F)\right) =0.
\]
Observe that, by Proposition 5.1.8 of \cite{hkst2015} (see Proposition \ref{prop:hkst}), an arc-length parametrized curve $\ga_s$ is necessarily absolutely continuous. Hence, for such
a curve, the notion of spending positive time in $F$ and the notion of having positive length in $F$ coincide.

We end this section by defining admissible densities for a family of curves in a metric space.
\begin{definition}
Let $\Ga$ be a family of locally rectifiable curves in the metric space $(X,d)$. A {\it density}, 
that is, a non-negative Borel function
$\rho: X\rightarrow[0,\infty]$ is {\it admissible} for $\Ga$ if
\[
\ell_\rho(\ga):=\int_\ga \rho ds\geq 1\qquad\forall \ga\in\Ga.
\]
\end{definition}
We write $\adm(\Ga)$ for the set of all admissible densities for $\Ga$, and note that this is purely a metric concept.

\subsection{The supremum-modulus on metric spaces}

On finite graphs, $\infty$-modulus is connected to shortest paths. Here, we extend this connection to the setting of metric spaces.

Given a family $\Ga$ of locally rectifiable curves in the metric space $(X,d)$, the {\it supremum-modulus} of $\Ga$ is 
\[
\Mod_{\rm sup}(\Ga):=\inf_{\rho\in\adm(\Ga)}\sup_X(\rho), 
\]
where $\sup_X(\rho)=\sup\{\rho(x): x\in X\}$. 
\begin{proposition}
Let $\Ga$ be a non-empty family of locally rectifiable curves in a metric space $(X,d)$.
Assume that
\[
0< \ell(\Ga):=\inf_{\ga\in\Ga} \length(\ga)<\infty.
\]
Then
\begin{equation}\label{eq:mod-sup-ell}
\Mod_{\rm sup}(\Ga) = \frac{1}{\ell(\Ga)}.
\end{equation}
\end{proposition}

\begin{pf}
Since  $\ell(\Ga)>0$, the density $\rho_0\equiv\frac{1}{\ell(\Ga)}$ is well-defined.
Note that for all $\ga\in\Ga$:
\[
\ell_{\rho_0}(\ga)=\int_\ga \rho_0 ds=\rho_0\length(\ga)=\frac{\length(\ga)}{\ell(\Ga)}\geq 1.
\]

Therefore, $\rho_0\in\adm(\Ga)$, hence $\adm(\Ga)\neq \emptyset$ and
\[
\Mod_{\rm sup}(\Ga)\leq \sup_{X}(\rho_0)=\ell(\Ga)^{-1}<\infty.
\]
Conversely, suppose $\rho\in\adm(\Ga)$. Then, given $\ga\in\Ga$,
\[
1\leq \int_\ga \rho ds = \int_0^{\length(\ga)}\rho(\ga_s(t))dt\leq\sup_X(\rho) \length(\ga).
\]
Since $\ga\in\Ga$ and $\rho\in\adm(\Ga)$ are arbitrary and both $\Ga$ and $\adm(\Ga)$ are non-empty, we can take the infimum  and get
\[
\ell(\Ga)\Mod_{\rm sup}(\Ga)\geq 1.
\]
\hfill\qedsymbol
\end{pf}
\begin{remark}
Note that (\ref{eq:mod-sup-ell}) makes sense also in limiting cases. For instance, if  
$\ell(\Ga)=\infty$, then any constant $\rho> 0$ is admissible, hence $\Mod_{\rm sup}(\Ga)=0$.
At the other extreme, if $\ell(\Ga)=0$, then there are arbitrarily short curves in $\Ga$, and hence for each
$\rho\in\adm(\Ga)$ we must have that $\sup_X(\rho)=\infty$ and hence $\Mod_{\rm sup}(\Ga)=\infty$.
\end{remark}

\section{Infinity modulus on metric measure spaces}\label{sec:infty-mod}

Geometric function theory grew out of complex analysis and real analysis on Euclidean spaces. It is therefore common to assume that the metric space $X$ is equipped with a regular Borel measure $\mu$.
The triple $(X,d,\mu)$ is referred to as a {\it metric measure space}.
For example, the measure $\mu$ can be obtained from the metric $d$ as a Hausdorff measure for an appropriate dimension, but it doesn't have to be. It is customary to require that $(X,d)$ is separable, and that for every point $x\in X$ there is a radius $r>0$ such the corresponding metric ball has  positive and finite measure, i.e.,  $0<\mu(B(x,r))<\infty.$ 

\subsection{Families of curves that are $\infty$-exceptional}

On a metric measure space $(X,d,\mu)$, it makes sense to talk about the {\it essential supremum}:
\[
\|\rho\|_\infty:=\inf\{a\geq 0: \mu(\{x: \rho(x)>a\})=0\}.
\]
Given a curve family $\Ga$, define
\[
\Mod_{\infty}^*(\Ga):=\inf_{\rho\in\adm(\Ga)}\|\rho\|_\infty. 
\]
Next, we establish some standard modulus properties for $\Mod_\infty^*$.
\begin{lemma}\label{lem:modprop}
$\Mod_{\infty}^*(\Ga)$ has the following properties:
\begin{itemize}
\item[(i)] If $\Ga_1\subset\Ga_2$, then  $\Mod_\infty^*(\Ga_1)\leq \Mod_\infty^*(\Ga_2)$ {\bf (monotone)}; 
\item[(ii)] $\Mod_\infty^*(\bigcup_{j\in\mathbb{N}}\Ga_j)\leq  
\sum_{j=1}^\infty\Mod_\infty^*(\Ga_j)$ {\bf (subadditive)}; 
\item[(iii)] If $\Ga_1\subset\Ga_2\subset\cdots$ and $\Ga=\bigcup_{j=1}^\infty \Ga_j,$
then $\lim_{j\rightarrow \infty}\Mod_\infty^*(\Ga_j)=\Mod_\infty^*(\Ga)$  {\bf (continuous from below)};
\item[(iv)] If for every $\ga\in\Ga_1$ there is a subcurve $\si\subset\ga$ such that $\si\in\Ga_2$ (denoted $\Ga_2\preceq\Ga_1$), then $\Mod_\infty^*(\Ga_1)\leq \Mod_\infty^*(\Ga_2)$ {\bf (shorter walks)}.
\end{itemize}
\end{lemma}
\begin{pf}
(i) Suppose that $\Ga_1\subset\Ga_2$. Then $\adm(\Ga_2)\subset \adm(\Ga_1)$. Thus, $\Mod_\infty^*(\Ga_1)\leq \Mod_\infty^*(\Ga_2)$.

(ii) Fix $\ep>0$, and suppose
that $\Ga=\bigcup_{j=1}^\infty \Ga_j$. For each $j$, find $\rho_j\in\adm(\Ga_j)$ such that
\[
\|\rho_j\|_\infty\leq \Mod_\infty^*(\Ga_j)+\frac{\ep}{2^j}.
\]
Let $\rho=\sum_j \rho_j$. Then $\rho\in\adm(\Ga)$ and 
\[
\Mod_\infty^*(\Ga)\leq \|\rho\|_\infty\leq \sum_{j=1}^\infty\|\rho_j\|_\infty \leq\sum_{j=1}^\infty\Mod_\infty^*(\Ga_j)+\ep
\]
Now let $\ep$ tend to zero.

(iii) By (i), the limit exists and is less than $\Mod_\infty^*(\Ga)$.  So, if the limit is infinite, we are done.  Otherwise, let $M:=\lim_{j\rightarrow\infty}\Mod_\infty^*(\Ga_j)$. For every $j=1,2,\dots$, there is $\rho_j\in\adm(\Ga_j)$ such that $\|\rho_j\|_\infty \leq M +1/j$. Define 
\begin{equation}\label{eq:suprho}
\rho^{(k)}(x):=\sup_{j\geq k} \rho_j(x).
\end{equation}
Also, for every $j=1,2,\dots$, there is a set $N_j$, with $\mu(N_j)=0$, such that $\rho_j< M+2/j$ outside of $N_j$. Then $N:=\cup_j N_j$ also has measure zero, and outside of $N$ we have $\rho^{(k)}<M+2/k$. In particular,
$\|\rho^{(k)}\|_\infty\leq M+2/k$.

Now fix $k=1,2,\dots$ and fix $\ga\in \Ga$. Then $\ga\in \Ga_j$ for some $j$. Let $i=\max\{j,k\}$.
Then since $i\geq j$, we have $\rho_i\in\adm(\Ga_j)$ and $\ell_{\rho_i}(\ga)\geq 1$.  Since $k\leq i$, by (\ref{eq:suprho}), $\ell_{\rho^{(k)}}(\ga)\geq 1$. So we have shown that $\rho^{(k)}\in\adm(\Ga)$, which implies that $\Mod_\infty^*(\Ga)\leq \|\rho^{(k)}\|_\infty\le M+2/k$. Letting $k$ tend to infinity, yields $\Mod_\infty^*(\Ga)\leq M$. 

The other direction, $\Mod_\infty^*(\Ga)\geq M$, follows from monotonicity (i). 

(iv) Note that $\adm(\Ga_2)\subset\adm(\Ga_1)$, so the infimum is taken over a larger set.
\hfill\qedsymbol
\end{pf}
\begin{definition}
\label{def:exceptional}
A curve family $\Ga$ is said to be {\rm $\infty$-exceptional}, if $\Mod_\infty^*(\Ga)=0$.
\end{definition}
\begin{example}\label{ex:loc-rect-exceptional}
Let $\Ga_\infty$ be the collection of all locally rectifiable curves that are not rectifiable. 
Then $\Ga_\infty$ is $\infty$-exceptional. Indeed, every constant density 
$\rho=a>0$ is admissible for $\Ga_\infty$. Therefore, $\Mod_\infty^*(\Ga_\infty)=0$.
\end{example}
As a consequence of Example \ref{ex:loc-rect-exceptional}, from now on we will always assume, without loss of generality, that our curve families consist uniquely of rectifiable curves.

In the next lemma we combine Lemma 5.7 and Lemma 5.8 of \cite{durandcartagena-jaramillo:jmaa2010}, and add the direction ${\rm (c)}\Rightarrow {\rm (d)}$. We rewrite the whole proof here for completeness. 
\begin{lemma}\label{lem:exceptional}
Let $\Ga$ be a family of rectifiable curves in the metric measure space $(X,d,\mu)$. Then, the following are equivalent:
\begin{itemize}
\item[(a)] $\Ga$ is $\infty$-exceptional.
\item[(b)] There exists $\rho: X\rightarrow [0,\infty)$ such that $\|\rho\|_\infty<\infty$ and $\ell_\rho(\ga)=\infty$ for all $\ga\in\Ga$.
\item[(c)] There exists $\rho: X\rightarrow [0,\infty)$ such that $\|\rho\|_\infty= 0$ and $\ell_\rho(\ga)=\infty$ for all $\ga\in\Ga$.
\item[(d)] There is a Borel set $F$ with $\mu(F)=0$ such that $\Ga\subset\Ga_F^\ell$. In words, there is a set of measure zero such that every curve in the family has positive length in that set.
\end{itemize}
\end{lemma}

\begin{pf}
${\rm (a)}\Rightarrow{\rm (b)}$: Assume $\Mod_\infty^*(\Ga)=0$. Then for $k=1,2,\dots$ there is $\rho_k\in \adm(\Ga)$ such that $\|\rho_k\|_\infty\leq 2^{-k}$. Set $\rho:=\sum_{k=1}^\infty\rho_k$. Then, $\|\rho\|_\infty\leq 1$ and $\ell_\rho(\ga)=\sum_{k=1}^\infty\ell_{\rho_k}(\gamma)=\infty$ for all $\ga\in \Ga$.

${\rm (b)}\Rightarrow{\rm (c)}$: We may assume without loss of generality that 
$\|\rho\|_\infty\leq 1$ and $\ell_\rho(\ga)=\infty$ for all $\ga\in\Ga$. 
For $a\geq 0$, consider the level sets
\[
S_a(\rho)=\{x: \rho(x)>a\}.
\] 
Since $\|\rho\|_\infty\leq 1$, we have $\mu(F)=0$, where $F=S_1(\rho)$.
Define $\tilde{\rho}:=\rho\  \ones_{F}$.
Then, $\|\tilde{\rho}\|_\infty=0$. We are left to show that $\ell_{\tilde{\rho}}(\ga)=\infty$ for all $\ga\in \Ga$. Using (b) and the rectifiability of $\ga$,
\begin{align*}
\infty & = \int_\ga \rho ds -\ell(\ga)= \int_{\ga\cap F}(\rho-1)ds+\int_{\ga\setminus F}(\rho-1)ds \\
&\le \int_{\ga\cap F}\rho ds =\int_\ga \tilde{\rho}ds.
\end{align*}

${\rm (c)}\Rightarrow{\rm (d)}$: Assume $\|\rho\|_\infty= 0$ and $\ell_\rho(\ga)=\infty$ for all $\ga\in\Ga$. Then $\mu(S_a(\rho))=0$ for every $a>0$.  Set $F=S_0(\rho)=\bigcup_{k\in\mathbb{N}}S_{1/k}(\rho)$. Then, by 
the subadditivity of measures,  
$\mu(F)=0$. However,  by the Cavalieri principle, for all $\ga\in\Ga$, 
\[
\infty=\int_\ga \rho ds= \int_0^{\length(\ga)} \rho(\ga_s(t))dt=\int_0^\infty m_1\left( \ga_s^{-1}(S_a(\rho))   \right)da.
\]
In particular, there is $a_0>0$, such that $m_1\left( \ga_s^{-1}(S_{a_0}(\rho))   \right)>0$.
Therefore, choosing a positive integer $k$ such that $a_0>1/k$, we have
\[
\int_\ga \ones_{\ga_s^{-1}(F)}\, du 
\ge \int_{0}^{\length(\ga)} \ones_{\ga_s^{-1}(S_{1/k}(\rho))}(t) dt>0.
\]
So $\Ga\subset\Ga_F^\ell$.

${\rm (d)}\Rightarrow{\rm (a)}$: Assume $\mu(F)=0$ and $\Ga\subset\Ga_F^\ell$.
For $k=1,2,\dots$, define 
\[
\Ga_k=\left\{\ga\in \Ga_F^\ell: \int_{\ga\cap F}ds\geq 1/k\right\}.
\]
Namely, these are the curves in $\Ga$ that have at least length $1/k$ in $F$.
Let $\rho_k:=k\ones_F$. Then $\rho_k\in\adm(\Ga_k)$ and $\|\rho_k\|_\infty=0$, so $\Mod_\infty^*(\Ga_k)=0$. Therefore, $\Mod_\infty^*(\Ga_F^\ell)=0$ by subadditivity (Lemma \ref{lem:modprop} (ii)) and 
$\Mod_\infty^*(\Ga)=0$ by monotonicity (Lemma \ref{lem:modprop} (i)).

Hence, we have shown that ${\rm (a), (b), (c), (d)}$ are all equivalent.
\hfill\qedsymbol
\end{pf}
\begin{remark}
If we allow $\rho$ to take on the value of $\infty$ at some points of $X$, then the proof of ${\rm (d)}\Rightarrow{\rm (a)}$ in Lemma \ref{lem:exceptional} is simplified: just set $\rho=\infty\, \ones_{F}$ and see that $\rho$ is admissible.
\end{remark}

\begin{definition}
We say that a property holds for {\rm $\infty$-almost every curve}, if it fails only for an $\infty$-exceptional set of curves.
\end{definition}
For instance, Lemma \ref{lem:exceptional} says that if $F$ is a Borel set with $\mu(F)=0$, then $\infty$-almost every curve has no length in $F$.

\subsection{Infinity modulus}

In addition to $\Mod_{\rm sup}$ and $\Mod_\infty^*$ here we consider a third notion of infinity-modulus.
We will show that the latter two notions are related.

\begin{definition}
We say that $\rho: X\rightarrow [0,\infty)$ is {\rm weakly admissible} and write $\rho\in\wadm(\Ga)$, if
\[
\int_\ga \rho ds\geq 1\qquad \text{for $\infty$-a.e. $\ga\in\Ga$}.
\]
Then, the {\rm $\infty$-modulus} of a family $\Ga$ is
\[
\Mod_\infty(\Ga):=\inf_{\rho\in\wadm(\Ga)}\|\rho\|_\infty.
\]
\end{definition}
\begin{remark}
Since it is easier to be weakly admissible than admissible,
\[
\Mod_\infty(\Ga)\leq \Mod_\infty^*(\Ga).
\]
In particular, an $\infty$-exceptional family has $\infty$-modulus zero. 
\end{remark}

\begin{lemma}
\label{lem:infty-mods-equal}
We have
\[
\Mod_\infty(\Ga) = \Mod_\infty^*(\Ga).
\]
\end{lemma}

\begin{pf}
In light of the above remark, it suffices to show that 
$\Mod_\infty(\Ga)\ge \Mod_\infty^*(\Ga)$. To this end, let $\rho$ be weakly admissible for $\Ga$. Then there is
a family $\Gamma_0\subset\Ga$ with $\Mod_\infty^*(\Gamma_0)=0$ such that whenever $\gamma\in\Ga\setminus\Gamma_0$
we have $\int_\gamma\rho\, ds\ge 1$. By Lemma~\ref{lem:exceptional} (c), there is a Borel function $\rho_0:X\to[0,\infty)$ such that
$\Vert\rho_0\Vert_\infty=0$ and such that for each $\gamma\in\Gamma_0$ we have $\int_\gamma\rho_0\, ds=\infty$. Note then that
$h:=\rho+\rho_0$ belongs to $\text{Adm}(\Ga)$. Thus 
\[
 \Mod_\infty^*(\Ga)\le \Vert h\Vert_\infty=\Vert \rho+\rho_0\Vert_\infty\le \Vert\rho\Vert_\infty+\Vert\rho_0\Vert_\infty=\Vert\rho\Vert_\infty.
 \]
 Taking the infimum over all $\rho$ that are weakly admissible for $\Ga$ yields that
\[
\Mod_\infty^*(\Ga)\le \Mod_\infty(\Ga),
\]
as desired.
\hfill\qedsymbol
\end{pf}

\subsection{Essential length}
\begin{definition}\label{def:ess-length}
Let $\Ga$ be a family of curves. For every $a\geq 0$, let
\begin{equation}\label{eq:levelset}
\Ga(a):=\{\ga\in\Ga: \ell(\ga)<a\}.
\end{equation}
The {\rm essential length} of $\Ga$ is
\[
{\rm ess}\ell(\Ga) := \sup\left\{a\geq 0: \text{$\Ga(a)$ is $\infty$-exceptional}\right\}.
\]
\end{definition}
Note that, by definition, we always have
\[
\ell(\Ga)\leq {\rm ess}\ell(\Ga).
\]
\begin{remark}\label{rem:essga}
For all $a<{\rm ess}\ell(\Ga)$, the family $\Ga(a)$ from (\ref{eq:levelset}) is $\infty$-exceptional. 
Writing $a_0:={\rm ess}\ell(\Ga)$, and using
the subadditivity of $\Mod_\infty^*$, Lemma \ref{lem:modprop} (ii), we get that $\Ga(a_0)$ is $\infty$-exceptional.
\end{remark}

\begin{theorem}\label{thm:mod-essl}
Let $(X,d,\mu)$ be a metric measure space with $\mu$ Borel. Let $\Ga$ be a family of  rectifiable curves in $X$.
\begin{itemize}
\item[(a)] If ${\rm ess}\ell(\Ga)\in(0,\infty)$,
then
\[
\Mod_\infty(\Ga)=\frac{1}{{\rm ess}\ell(\Ga)}\in (0,\infty),
\]
and $\rho_0\equiv {\rm ess}\ell(\Ga)^{-1}$ is an extremal weakly admissible density.

\item[(b)] If ${\rm ess}\ell(\Ga)=0$, then $\Mod_\infty(\Ga)=\infty$ and no extremal weakly admissible density exists in $L^\infty(X)$.

\item[(c)] If ${\rm ess}\ell(\Ga)=\infty$, then $\Mod_\infty(\Ga)=0$.
\end{itemize}
\end{theorem}

\begin{pf}
Assume that $a_0:={\rm ess}\ell(\Ga)\in(0,\infty)$. Suppose that
$\rho\in\wadm(\Ga)$. Set $F=\{x: \rho(x)>\|\rho\|_\infty\}$. Then $F$ is a Borel set with $\mu(F)=0$, so, by Lemma \ref{lem:exceptional},  
$\Ga_F^\ell$ 
is $\infty$-exceptional.
Also, by definition of weakly admissible, the family 
\[
\tilde{\Ga}=\left\{\ga\in\Ga: \int_\ga \rho ds< 1\right\} 
\]
is $\infty$-exceptional.
Finally,   as seen in Remark \ref{rem:essga}, $\Ga(a_0)$ is $\infty$-exceptional.
Therefore, 
by subadditivity $\Ga_F^\ell\cup \tilde{\Ga}\cup \Ga(a_0)$ is $\infty$-exceptional. On the other hand, 
${\rm ess}\ell(\Ga)<\infty$ implies that $\Ga$ is not $\infty$-exceptional, by Definition \ref{def:ess-length}. So we can find at least one curve $\ga\in\Ga\setminus(\Ga_F^\ell\cup \tilde{\Ga}\cup \Ga(a_0))$. For these curves, we have
\[
1\leq \int_\ga \rho ds =\int_{\ga\cap F}\rho ds+\int_{\ga\setminus F}\rho ds=\int_{\ga\setminus F}\rho ds\leq \|\rho\|_\infty\ell(\ga).
\]
Moreover,  for every $\ep>0$, $\Ga(a_0+\ep)$ is also not $\infty$-exceptional, hence repeating the argument above we can also find a curve $\ga\in\Ga(a_0+\ep)$ such that
\[
1\le \|\rho\|_\infty\ell(\ga)\le \|\rho\|_\infty({\rm ess}\ell(\Ga)+\ep).
\]
Letting $\ep\rightarrow 0$, we find that
\[
1\leq \|\rho\|_\infty {\rm ess}\ell(\Ga).
\]
Since $\rho\in \wadm(\Ga)$ was arbitrary, we obtain that
\[
\Mod_\infty(\Ga)\geq\frac{1}{{\rm ess}\ell(\Ga)}.
\]

Conversely, let $\rho_0\equiv {\rm ess}\ell(\Ga)^{-1}$. Then,  for every $\ga\in \Ga\setminus\Ga(a_0)$, we have 
$\ell(\ga)\geq {\rm ess}\ell(\Ga)$. Thus
\[
\int_\ga \rho_0 ds =\frac{\ell(\ga)}{{\rm ess}\ell(\Ga)}\geq 1.
\]
By Remark \ref{rem:essga}, $\Ga(a_0)$ is $\infty$-exceptional. So $\rho_0$ is weakly admissible and therefore,
\[
\Mod_\infty(\Ga)\leq\|\rho_0\|_\infty=\frac{1}{{\rm ess}\ell(\Ga)}.
\]
Thus (a) is proved.

For (b), assume that ${\rm ess}\ell(\Ga)=0$.  We shall show that no
weakly admissible density $\rho$ has finite essential norm; by
convention, the infimum of an empty set is $\infty$, so it will follow
that $\Mod_\infty(\Gamma)=\infty$.  Let $\rho$ be a density such that
$\|\rho\|_\infty<\infty$ and, as before, define
$F:=\{z:\rho(z)>\|\rho\|_\infty\}$.  By
Definitions~\ref{def:exceptional} and~\ref{def:ess-length} and
Lemma~\ref{lem:infty-mods-equal}, ${\rm ess}\ell(\Ga)=0$ implies that
$\Mod_\infty(\Ga(a))>0$ for any $a>0$.  Moreover, since $\Ga^\ell_F$
is $\infty$-exceptional, it follows that
$\Mod_\infty(\Ga(a)\setminus\Ga^\ell_F)>0$.  But, for any
$\gamma\in\Ga(a)\setminus\Ga^\ell_F$,
\begin{equation*}
  \int_\ga \rho\,ds \le \|\rho\|_\infty\ell(\gamma) < a\|\rho\|_\infty.
\end{equation*}
For sufficiently small $a$,
\begin{equation*}
  \Ga(a)\setminus\Ga^\ell_F \subset \left\{\gamma\in\Gamma : \int_\gamma\rho\;ds < 1\right\}.
\end{equation*}
Since the former set has positive $\infty$-modulus, $\rho\notin\wadm(\Ga)$.

Finally for (c), assume that ${\rm ess}\ell(\Ga)=\infty$. Then, $\Ga(n)$ is $\infty$-exceptional for every $n\in\N$. Also,
 $\Ga=\bigcup_{n\in\mathbb{N}}\Ga(n)$. Therefore, by subadditivity of modulus, we have that
$\Mod_\infty(\Ga)=0$.

\hfill\qedsymbol
\end{pf}

\section{The essential metric}\label{sec:ess-metric}
Consider the ``connecting'' families $\Ga(x,y)$ consisting of all curves connecting two points $x\neq y$. Define
\begin{equation}\label{eq:ess-metric}
d_{\rm ess}(x,y):=\Mod_\infty(\Ga(x,y))^{-1}={\rm ess}\ell(\Ga(x,y)).
\end{equation}
Note that $d_{\rm ess}(x,y)$ could be infinite for some $x,y\in X$, e.g., if $\Ga(x,y)=\emptyset$. 
\begin{example}
 If $X$ is the Sierpinski gasket in the
plane, equipped with the Euclidean metric and the natural Hausdorff measure, then 
from the results of~\cite{Durand-Sh-Williams} (see also~\cite{Bourdon-Pajot}),
the collection of all
non-constant rectifiable curves in $X$ is $\infty$-exceptional; thus in this case as well, even though $\Ga(x,y)$ is non-empty for each pair of points
$x,y\in X$, we have that $d_{\rm ess}(x,y)=\infty$ when $x\ne y$.
\end{example}
\begin{theorem}\label{thm:ess-metric}
The function $d_{\rm ess}: X\times X\rightarrow \R_{\ge 0}\cup\{\infty\}$ defined in (\ref{eq:ess-metric}) for $x\neq y$, and defined to be zero on the diagonal of $X\times X$, is an extended metric on $X$. 
Moreover, if for each $x,y\in X$ with
$x\ne y$ we have $\Mod_\infty(\Ga(x,y))>0$, then $d_{\rm ess}$ is a metric on $X$, with the property that $d_{\rm ess}\ge d$.
\end{theorem}
\begin{pf}
Let $x,y\in X$ with $x\ne y$. 
Since every curve can be reversed, the symmetry axiom holds: $d_{\rm ess}(x,y)=d_{\rm ess}(y,x)$.

Next we show that if $x\ne y$, then $d_{\rm ess}(x,y)\ge d(x,y)>0$. Since $X$ is a metric space, we have $d(x,y)>0$. By definition of length, every curve $\ga$ connecting $x$ to $y$ satisfies
$\ell(\ga)\ge d(x,y)>0$. Therefore, for any $a<d(x,y)$, the family $\Ga(x,y)(a)$, using the notation from~\eqref{eq:levelset}, is empty, and thus $\infty$-exceptional. This shows that 
$d_{\rm ess}(x,y)\ge d(x,y)>0$. 

Next, we verify the triangle inequality. 
Fix three distinct points $a,b,c\in X$. 
For simplicity, let $\Ga_1=\Ga(a,c)$, $\Ga_2=\Ga(c,b)$ and $\Ga_0=\Ga(a,b)$. 
Also assume that $\de_0:=d_{\rm ess}(a,b),\ \de_1:=d_{\rm ess}(a,c),\ \de_2:=d_{\rm ess}(c,b)$, are all positive.
Let $\Ga=\Ga(a,b;c)$ be the family of all curves that start at $a$, end at $b$, and go through $c$.  Then $\Ga\subset\Ga_0$, and so $\de_0={\rm ess}\ell(\Ga_0)\le {\rm ess}\ell(\Ga)$.

If $\de_1+\de_2=\infty$, then clearly $\de_0\le \de_1+\de_2$. So we can assume that both $\de_1$ and $\de_2$ are finite. Fix $\ep>0$ and let $\la:=\de_1+\de_2+2\ep$. 
We want to show that 
\begin{equation}\label{eq:modga-la}
\Mod_\infty(\Ga(\la))>0,
\end{equation} 
because, by monotonicity, that implies that $\Mod_\infty(\Ga_0(\la))>0$, and thus, by definition of essential length, $\de_0\le\la=\de_1+\de_2+\ep$. Then letting $\ep$ tend to zero, yields the conclusion.

To that end, fix a density $\rho$ which is  admissible 
for $\Ga(\la)$. Such a $\rho$ exists, because we have assumed that $\de_0>0$, so $\Mod^*_\infty(\Ga_0)<\infty$.
Fix two arbitrary curves:  $\alpha\in\Ga_1(\de_1+\ep/2)$ and 
$\beta\in\Ga_2(\de_2+\ep/2)$. Such curves always exist because, by Remark \ref{rem:essga},
\begin{equation}\label{eq:modgaj-pos}
\Mod_\infty(\Ga_j(\de_j+\ep/2))>0,\qquad \text{for $j=1,2$.}
\end{equation}
Write $\ga:=\alpha+\beta$ for the concatenation of $\alpha$ and $\beta$.
Then 
\[
\ell(\ga)=\ell(\alpha)+\ell(\beta)<(\de_1+\ep/2)+(\de_2+\ep/2).
\]
So $\ga\in\Ga(\lambda)$. Therefore, by admissibility,
$\int_{\ga}\rho\, ds\ge 1$. 
In particular,  if, for $j=1,2$, we set
\[
\ell_\rho(\Ga_j(\de_j+\ep/2)):=\inf_{\ga'\in\Ga_j(\de_j+\ep/2)}\int_{\ga'}\rho\, ds,
\] 
then
\begin{equation}\label{eq:lb}
\ell_\rho(\Ga_1(\de_1+\ep/2))+\ell_\rho(\Ga_2(\de_2+\ep/2))\geq 1.
\end{equation}

On the other hand, since such curves $\alpha$ and $\beta$ exist, $\ell_\rho(\Ga_j(\de_j+\ep/2))<\infty$,  for $j=1,2$. Moreover, we claim that 
\begin{equation}\label{eq:ell-rho-upperbound}
 \ell_\rho(\Ga_j(\de_j+\ep/2))
 \le \frac{\Vert\rho\Vert_\infty}{\Mod_\infty^*(\Ga_j(\de_j+\ep/2))},\qquad\text{for $j=1,2$}.
\end{equation}
To prove this, assume first that $\ell_\rho(\Ga_j(\de_j+\ep/2))>0$, for both $j=1,2$.
Then,
\begin{equation}\label{eq:rho-adm-j}
\frac{\rho}{\ell_\rho(\Ga_j(\de_j+\ep/2))}\in\adm(\Ga_j(\de_j+\ep/2)),\qquad\text{for $j=1,2$}.
\end{equation}
 In particular,
\[
0<\Mod_\infty^*(\Ga_j(\de_j+\ep/2))
\le \frac{\Vert\rho\Vert_\infty}{\ell_\rho(\Ga_j(\de_j+\ep/2))},\qquad\text{for $j=1,2$},
\]
which implies (\ref{eq:ell-rho-upperbound}). Lastly, if, say
$\ell_\rho(\Ga_1(\de_1+\ep/2))=0$, then (\ref{eq:ell-rho-upperbound}) holds trivially for $j=1$, and  (\ref{eq:lb}) implies that $\ell_\rho(\Ga_2(\de_2+\ep/2))\ge 1$, so the same admissibility argument in (\ref{eq:rho-adm-j}) works for $j=2$.

Combining (\ref{eq:lb}) and (\ref{eq:ell-rho-upperbound}), we get that 
\[
1\le \Vert\rho\Vert_\infty\left[\frac{1}{\Mod_\infty^*(\Ga_1(\de_1+\ep/2))}
  +\frac{1}{\Mod_\infty^*(\Ga_2(\de_2+\ep/2))}\right],
\]
and taking the infimum over all such admissible $\rho\in\adm(\Ga(\la))$ yields
\[
\frac{1}{\Mod_\infty^*(\Ga(\lambda))}\le 
   \bigg[\frac{1}{\Mod_\infty^*(\Ga_1(\de_1+\ep/2))}
  +\frac{1}{\Mod_\infty^*(\Ga_2(\de_2+\ep/2))}\bigg]<\infty,
\]
where the finiteness follows from (\ref{eq:modgaj-pos}).
This shows that $\Mod_\infty^*(\Ga(\lambda))>0$, hence (\ref{eq:modga-la}) is established,
and the triangle inequality is proved.

Finally, suppose that for each $x,y\in X$ with
$x\ne y$ we have $\Mod_\infty^*(\Ga(x,y))>0$. Then to show that $d_{\rm ess}$ is a
metric on $X$ it is enough to show that $d_{\rm ess}(x,y)<\infty$ for $x,y\in X$.
Note that by Lemma \ref{lem:modprop} (iii), $\Mod_\infty^*(\Ga(x,y))=\lim_{n\to\infty}\Mod_\infty^*(\Ga(x,y)(n))$. Therefore, since $\Mod_\infty^*(\Ga(x,y))>0$, we must have 
$\Mod_\infty^*(\Ga(x,y)(n))>0$ for some $n\in\N$. This implies that
$d_{\rm ess}(x,y)\le n<\infty$, completing the proof.
\hfill\qedsymbol
\end{pf}

The following definition is 
from~\cite{Durandcartagena-Jaramillo-Shanmugalingam:CVGMT2016}, and is
due to De Cecco and Palmieri~\cite{DeCecco-Palmieri}.

\begin{definition}
Let $(X,d,\mu)$ be a metric measure space.
Given a set $N\subset X$ with $\mu(N)=0$, and for $x,y\in X$ with $x\ne y$, we set
\[
d_N(x,y):=\inf\{\ell(\gamma)\, :\, \gamma\text{ connects }x\text{ to }y, \text{ and }
  m_1(\gamma_s^{-1}(N))=0\}.
\]
Also, define $\widehat{d}:X\times X\to\R$ by $\widehat{d}(x,x)=0$ and for $x\ne y$,
\[
\widehat{d}(x,y):=\sup\{d_N(x,y)\, :\, N\subset X\text{ with }\mu(N)=0\}.
\]
\end{definition}

As in the case of the essential metric $d_{\text ess}$, it might happen that $\widehat{d}$ is not finite. It was shown in~\cite{Durandcartagena-Jaramillo-Shanmugalingam:CVGMT2016}
that if $\mu$ is doubling and $X$ is complete, then 
$\widehat{d}$ is biLipschitz equivalent to the original
metric $d$ if and only if $X$ supports an $\infty$-Poincar\'e inequality.

\begin{remark}
From the results in~\cite{Durandcartagena-Jaramillo-Shanmugalingam:CVGMT2016}
it follows that if $\widehat{d}$ is a metric on $X$ such that $\mu$ is doubling with respect to this metric, then
$(X,\widehat{d},\mu)$ supports an $\infty$-Poincar\'e inequality. Here we say that $\mu$ is \emph{doubling} if there is 
a constant $C\ge 1$ such that whenever $x\in X$ and $r>0$, we have
\[
\mu(B(x,2r))\le C\, \mu(B(x,r)).
\]
We say that $X$ supports \emph{an $\infty$-Poincar\'e inequality}
if there are constants $C>0$ and $\lambda\ge 1$ such that whenever $u\in L^\infty(X)$ and $g$ is an upper gradient of $u$
\[
\frac{1}{\mu(B(x,r))}\int_{B(x,r)}|u-u_{B(x,r)}|\, d\mu\le C\, r\, \|g\chi_{B(x,\lambda r)}\|_\infty\qquad\forall x\in X, r>0.
\]
A non-negative Borel function $g$ on $X$ is said to be an \emph{upper
gradient} if for each rectifiable curve $\ga$ in $X$ we have
\begin{equation}\label{eq:upper-gr}
|u(x)-u(y)|\le \int_\ga g\, ds.
\end{equation}
The notion of upper gradients
is due to Heinonen and Koskela, see~\cite{hkst2015}. See the papers~\cite{Durandcartagena-Jaramillo-Shanmugalingam:CVGMT2016, Durand-Sh-Williams} for more on the $\infty$-Poincar\'e inequality. Heuristically, the $\infty$-Poincar\'e inequality gives us a way
of controlling the variance of a function on a ball in terms of its $\infty$-energy  on a slightly enlarged ball.
\end{remark}

\begin{proposition}
If $(X,d,\mu)$ is a metric measure space, then $\widehat{d}=d_{\text ess}$.
\end{proposition}

\begin{pf}
Fix $x,y\in X$ with $x\ne y$. We will first show that 
$\widehat{d}(x,y)\le \dess(x,y)$. To this end, note that if
$\dess(x,y)=\infty$, then the above inequality holds trivially. Therefore,
let us assume that $\dess(x,y)<\infty$. Recall that $\Ga:=\Ga(x,y)$ denotes the collection
of all rectifiable curves in $X$ connecting $x$ to $y$.
For each $\varepsilon>0$, if $a=[1+\varepsilon]\dess(x,y)$, then $\Mod_\infty(\Ga(a))>0$.
So, by Lemma \ref{lem:exceptional}(d), whenever $N\subset X$ with $\mu(N)=0$ there must exist a curve $\ga$ in
$\Ga(a)$ that does not spend positive time in $N$.
In particular,  $m_1(\ga_s^{-1}(N))=0$.  Hence,
\[
d_N(x,y)\le \ell(\ga)< [1+\varepsilon]\dess(x,y).
\]
Taking the supremum over all 
such nulls sets $N$ gives 
$\widehat{d}(x,y)\le [1+\varepsilon]\dess(x,y)$. 
Letting $\varepsilon\to 0$ gives the desired inequality.

We next show that $\dess(x,y)\le \widehat{d}(x,y)$. 
First suppose that 
$\dess(x,y)=\infty$. Then, by Theorem \ref{thm:mod-essl} (c),
$\Mod_\infty(\Ga(x,y))=0$. Thus by Lemma~\ref{lem:exceptional}(d), there is some $N\subset X$ such that
$\mu(N)=0$ and $\Ga(x,y)\subset\Ga_N^\ell$. Thus $\widehat{d}(x,y)\ge d_N(x,y)=\infty$.

Now consider the case that $\dess(x,y)<\infty$. By Theorem \ref{thm:ess-metric}, $\dess(x,y)\ge d(x,y)>0$, since $x\ne y$. Therefore,
$\Mod_\infty(\Ga(x,y)(a))=0$, whenever $a=[1-\varepsilon]\dess(x,y)$ for $0<\varepsilon<1$.
By Lemma~\ref{lem:exceptional}(d), there is a Borel set $N\subset X$, with $\mu(N)=0$, such that $\Ga(x,y)(a)$ has positive length in $N$. In particular, $d_N(x,y)\ge a$. Hence,
\[
\widehat{d}(x,y)\ge d_N(x,y)\ge [1-\varepsilon]\dess(x,y).
\]
Letting $\varepsilon\to 0$ gives $\widehat{d}(x,y)\ge \dess(x,y)$ as
desired.
\hfill\qedsymbol
\end{pf}
\begin{remark}
By Theorem \ref{thm:ess-metric}, $\dess\ge d$. Thus, if $\{x_i\}_i$ is a sequence in $X$ converging
with respect to $\dess$ to a point $x\in X$, then this sequence also converges in the original metric $d$. Therefore, in general,
the topology generated by $\dess$ is finer than that of $d$.

Given a metric measure space, it is of interest to know whether there are a great many curves
of controlled length connecting any given pair of disjoint continua. The quantity $\dess$ is a valuable tool in this, as
$\dess(x,y)=\infty$ tells us that there are very few such curves connecting $x$ to $y$. Therefore $\dess$ is a more sensitive
metric on $X$, and it would be interesting to know when the space $(X,\dess,\mu)$ supports an $\infty$-Poincar\'e inequality.
\end{remark}

\bibliographystyle{acm}
\bibliography{InfinityMod}
\def\cprime{$'$}

\end{document}